\documentclass[conference]{IEEEtran}
\IEEEoverridecommandlockouts
\usepackage{cite}
\usepackage{amsmath,amssymb,amsfonts}
\usepackage{amsmath}
\usepackage{algorithmic}
\usepackage{graphicx}
\usepackage{textcomp}
\usepackage{bm} 
\usepackage{xcolor}
\usepackage{booktabs}
\usepackage{subfig}

\def\BibTeX{{\rm B\kern-.05em{\sc i\kern-.025em b}\kern-.08em
T\kern-.1667em\lower.7ex\hbox{E}\kern-.125emX}}
\begin{document}

\title{A Traffic Evacuation Model for Enhancing Resilience During Railway Disruption}




\author{
    \IEEEauthorblockN{
        Hangli Ge\IEEEauthorrefmark{1}, Xiaojie Yang\IEEEauthorrefmark{2}, Jinyu Chen\IEEEauthorrefmark{1}, Francesco Flammini\IEEEauthorrefmark{3,4}, Noboru Koshizuka\IEEEauthorrefmark{1}
    }
    \IEEEauthorblockA{
        \IEEEauthorrefmark{1}Interfaculty Initiative in Information Studies, The University of Tokyo, Tokyo, Japan\\
        \IEEEauthorrefmark{2}Graduate School of Interdisciplinary Information Studies, The University of Tokyo, Tokyo, Japan\\
        \IEEEauthorrefmark{3}Department of Mathematics and Computer Science Ulisse Dini, University of Florence, Florence, Italy,\\
        and IDSIA USI-SUPSI, University of Applied Sciences and Arts of Southern Switzerland, Lugano, Switzerland\\
        Email: \{hanglige,xiaojieyang\}@g.ecc.u-tokyo.ac.jp,\\
        miraclec@csis.u-tokyo.ac.jp, francesco.flammini@unifi.it,
        noboru@koshizuka-lab.org
    }
}
\maketitle

\begin{abstract}
This paper introduces a traffic evacuation model for railway disruptions to improve resilience. The research focuses on the problem of failure of several nodes or lines on the railway network topology. We proposed a holistic approach that integrates lines of various operator companies as well as external geographical features of the railway system. The optimized evacuation model was mathematically derived based on matrix computation using nonlinear programming. The model also takes into account the capacity of the surrounding evacuation stations, as well as the travel cost. Moreover, our model can flexibly simulate disruptions at multiple stations or any number of stations and lines, enhancing its applicability. We collected the large-scale railway network of the Greater Tokyo area for experimentation and evaluation.  We simulated evacuation plans for several major stations, including Tokyo, Shinjuku, and Shibuya. The results indicate that the evacuation passenger flow (EPF) and the average travel time (ATT) during emergencies were optimized, staying within both the capacity limits of the targeted neighboring stations and the disruption recovery time.
\end{abstract}





\section{Introduction}
Resilience is a concept widely applied across fields such as economics, psychology, and engineering. It originally derived from the Latin word “resiliere”, was first defined in the ecological system \cite{holling1973resilience}.
Transportation systems, being at the forefront of various environmental and socioeconomic shocks, have made the study of resilience a topic of significant interest among transportation scholars. In particular, railway systems, as critical infrastructure, play a fundamental role in passenger or cargo transport, significantly contributing to economic growth, urban development, and environmental sustainability. 


However, under abnormal or emergency conditions—such as natural disasters, equipment failures, cyber attacks, operational disruptions, or accidents—trains may fail to arrive or depart on schedule or even experience service interruptions, leading to significant variations in railway operations. These issues can result in congestion, stoppages, or accidents, ultimately reducing passenger service quality and causing economic and productivity losses. Therefore, enhancing the resilience of transportation systems is essential for mitigating the impact of such disruptions
\cite{tang2022literature,zhou2019resilience}. 

\begin{figure}[t]
\centering
\includegraphics[height=2.1in, width=0.85\linewidth]{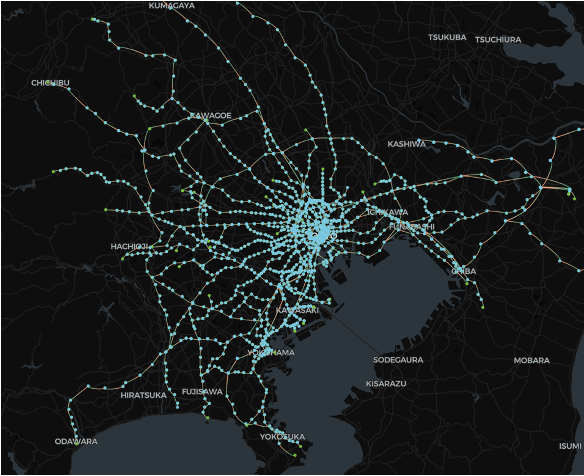}
\caption{The railway network of Greater Tokyo area }
\label{fig:net}
\end{figure}
Improving the resilience of railways to recover to normal conditions more quickly is a critical problem in the management of railway systems. Among various recovery strategies, traffic analysis and evacuation simulation are particularly significant. However, these tasks are highly complex due to the lack of readily available large-scale data on railway networks and passenger flows—whether at the cross-city, nationwide, or cross-country level. If the model focuses solely on a single railway line or lines operated by the same company, the process is relatively straightforward. However, constructing a comprehensive rail network presents significant challenges. The involvement of multiple railway operators further complicates data collection, integration, computation, and the simulation of decision-making within a large-scale railway system.

In this research, we propose a traffic evacuation model to enhance railway resilience during disruptions. \textbf{The model addresses the failure of multiple nodes or lines within the network topology. We introduce an optimization model for traffic evacuation that considers delay times, overcrowded passenger flows, and the capacities of neighboring stations.} We constructed a large-scale railway network dataset for the Greater Tokyo area (visualized in Fig.~\ref{fig:net}), integrating data from multiple operating companies and incorporating geographical features such as connectivity and distance, which were modeled in the evacuation cost matrix. Our traffic evacuation model employs an optimization algorithm based on matrix calculations, ensuring both intuitiveness and interpretability. The key contributions of this research can be summarized as follows:


\begin{itemize}
    \item \textbf{Proposal of a fusion approach} which integrates multiple geographic features to construct an evacuation cost matrix that dynamically adjusts based on predefined disruption durations.
    \item \textbf{Development of an optimization model} using nonlinear programming that accounts for the capacity of surrounding stations and travel costs. Additionally, the proposed model can flexibly simulate disruptions affecting multiple stations or spanning cross multiple lines, enhancing its applicability. Built on a solid mathematical foundation, the model leverages straightforward matrix calculations and optimization solutions, making it intuitive, easy to understand, and simple to deploy.
   \item \textbf{Validation through a case study} on Greater Tokyo’s large-scale railway network. Geographic features and passenger statistics were incorporated into simulations for several major stations. The optimization focused on evacuation passenger flow (EPF) and the average travel time (ATT), ensuring both remained within the capacity limits of the evacuated target stations and adhered to the designated disruption time.
\end{itemize}




\section{Related Work}

\subsection{Resillience Assessment}
A resilient transport system is crucial for maintaining functionality and accessibility, especially in the face of unexpected events, such as natural disasters, accidents or infrastructure failures \cite{chen2024resilience}. However, the operationalization of resilience practices in transport planning and found that resilience is not formally considered \cite{esmalian2022operationalizing}. Researchers have developed various models and systems to address railway resillience issues and continue to explore new methods for improving passenger experience and reducing delay or congestion when abnormaly events occur on the railway system. However, regarding the resillience, little is known about the performance dynamics during disruptions\cite{knoester2024data}. 

Data-driven approach to quantify the resilience of the railway network through assessment, utilizing historical traffic data from the Netherlands to reconstruct resilience curves and analyze disruption dynamics using composite indicators and statistical tests, providing valuable insights for improving system resilience\cite{knoester2024data}. Similar research has evaluated the resilience of complementary networks in various disaster scenarios from the perspective of network structure and function \cite{chen2024resilience}. These resilience-related studies are constructive, but one limitation in resilience measurement is that the features of the network have not been fully considered.

\subsection{Practice for Improving Resillience}
Several practices are used to improve the resilience of railway systems. They put the focus on transport planning mainly include the identification of critical nodes and edges \cite{lee2022quantitative}, transport network design \cite{xu2021enhancing}, emergency response, and evacuation \cite{sohouenou2021assessing,esmalian2022operationalizing}. Digital twins can help optimize maintenance, reduce downtime, and improve safety by providing real-time monitoring and analysis of system performance \cite{Kaewunruen2022}. In addition,  the use of sensors and data analytics or AI (Artificial intelligence)-based predictions can help predict when maintenance is required, reducing the likelihood of equipment failure and improving overall system reliability \cite{bevsinovic2021artificial}.  

Similarly, a study proposes a novel method to optimize the resilience of urban rail systems, focusing on the combined impact of delayed trains and overcrowded passenger flows \cite{li2023optimization}. Additionally, the use of complex network principles can help analyze the network structure and identify potential bottlenecks, enabling targeted improvements. ``KUTTY" simulator, which simulates the flow of passengers quantitatively in all segments of the network \cite{Tomoeda2012}.  The model has also been used to simulate the effects of a virtual accident at Otemachi Station, showing that the flow pattern of passengers will change significantly\cite{Tomoeda2012}. 


Traffic evacuation models also have been proposed for optimizing the evacuation process, minimize travel delays, and ensure the safety of passengers. An approach is to use a resilience-based optimization model, which combines resilience evaluation and dispatching optimization to maximize the resilience index of evacuated passengers of bus system \cite{Zhang2023}. A similar approach is to represent an integrated network, generate candidate bus-bridging routes using the K-shortest paths algorithm, and solve the optimization model to determine the optimal vehicle allocation among the candidate routes \cite{Jaber2024}. The proposal was applied to a case study in the Ile de France region, Paris and suburbs, to optimize transportation management during interruptions.

In summary, data-driven approaches for simulating traffic evacuation are significant to improve the resilience of railway systems. However, it requires a holistic view that considers various components, including infrastructure, external features of users or cities, operator organization, integration level, for emergency risk management and recovery. This needs to be achieved through the development of integrated data management systems that bring together various stakeholders, including operators, maintenance personnel such as local government, and IT-service providers and so on. More features of the rail network are modeled, the more the constraints of the optimization solutions align with real-world conditions, making the 
optimization solutions more feasible.

\section{Problem Definition and Notations}


We propose a traffic evacuation model to improve resilience during rail disruption. The railway traffic evacuation problem was mathematically defined and optimized using nonlinear programming. 

\subsection{Input of Railway Network}
Suppose we have $n$ stations, represented by the set of station nodes of \{$v_{1},v_{2},v_{n-1},v_{n}$\}. Hereafter, we define the following notations:
\[
\begin{aligned}
\bm{Y},\, \bm{J},\, \bm{\hat{J}},\, \bm{X} \in \mathbb{R}^{n \times 1}, 
\quad
\bm{A_{con}},\bm{A_{dis}},\bm{A_{cost}} \in \mathbb{R}^{n \times n},\bm{T_{lm}}>0\\
\end{aligned}
\]
$\bm{Y}$, in which each element \( y_i \), represents the passenger count at each station. $\bm{J}$ is an all-ones matrix, and $\bm{\hat{J}}$ is used to indicate whether a station is normal or blocked, with a value of 1 for blocked stations and 0 for normal stations. $\bm{X}$, where each element \( x_i \), represents the multiplicative increase in station volume capacity (i.e., a value greater than 1). $\bm{T_{lm}}$ denotes the disrupted duration (in minutes), which is greater than 0.

\subsection{Dynamic Cost Matrix Generation}

The geographical features of the railway network can be defined as the edge values in a graph, which can be represented as elements of matrix as follows. 

\textbf{Connection Matrix of $\bm{A_{con} \in \mathbb{R}^{n \times n}}$ }: Connectivity is the most fundamental principle that ensures the completeness and correctness of a railway network. Each element $a_{i,j}$ in $\bm{A_{con}}$ is denoted by Equation~\ref{eq:w_c}.\begin{equation}
a_{i,j} =
\label{eq:w_c} 
\begin{cases}
     1   &   \text{if $v_{i}$ connects to $v_{j}$ in the railway line}\\  
     0  & otherwise\\
\end{cases}
\end{equation}



\textbf{Distance Matrix of $\bm{A_{dis}} \in \mathbb{R}^{n \times n}$}: Let $\bm{A_{dis}}$ represents the distance matrix, where each element $a_{i,j}$ in $\bm{A_{dis}}$ is given in Equation ~\ref{eq:d_lat_lon}. Here, $d_{i,j}$ denotes the Haversine distance between two points on the globe, identified by their latitude and longitude. R is the Earth’s radius, approximately 6371 km.
\begin{equation}
 \begin{aligned}
a_{i,j} = d_{i,j}= 2\bm{R}*\textrm{argsin}(sqrt((\sin^{2}\frac{(lat_{i}-lat_{j})}{2}\\+\cos(lat_{i})\cos(lat_{j})\sin^{2}\frac{(lon_{i}-lon_{j})}{2})))
 \end{aligned}
\label{eq:d_lat_lon} 
\end{equation}


\textbf{Fused Travel Cost Matrix of $A_{cost}$}: In this research, the travel cost matrix was constructed based on two types of basic matrix: connection matrix and distance matrix to exploit more different but useful heterogeneous spatial correlation among stations. The fused cost matrix (denoted as $A_{cost}$) were calculated by fusing the $A_{con}$  and $A_{dis}$ matrices at the element level, as shown in Equation~\ref{eq:a_cost}. 
\begin{equation}
\begin{aligned}
    \bm{A_{cost}} &= \bm{T_{train}} \odot \bm{A}_{con} \;+\; (\bm{I} \;-\; \bm{A}_{con}) \;\odot\; \left(\frac{\bm{A}_{dis}}{s}\right) \\
    &\quad \text{with} \quad
    \bm{A_{cost}} = \begin{cases} 
      \infty, & \text{if } a_{ij} > T_{lm}, \\
      a_{ij}, & \text{otherwise},
    \end{cases}
\end{aligned}
\label{eq:a_cost}
\end{equation}

Throughout the formulation, the operator ``$\odot$'' denotes the Hadamard (element-wise) product. $\bm{I}$ is the identity matrix. $\bm{T_{train}} \in \mathbb{R}^{n \times n}$ where each element $t_{i,j}$ denotes the shortest travel time by train from station $v_{i}$ to $v_{j}$. To determine the travel cost, we calculate the minimum travel time between all pairs of stations with one transfer (either walking or taking a train). $\left(\frac{\bm{A}_{dis}}{s}\right)$ denotes the walking time between two stations and we set the average walking speed is $s$ (set to be 5km~/h). In particular, during our optimization, we consider it impractical for evacuation time to exceed the station's recovery time. Therefore, we set the entry $a_{i,j}$ in $\bm{A_{cost}}$ that exceeds a certain disruption duration of $T_{lm}$ (e.g., 30 minutes or 60 minutes) to infinity, thus helping us quickly narrow down to reasonable candidate stations during the computation.


\begin{figure*}[t]
\centering
\includegraphics[height=2.1in, width=1\linewidth]{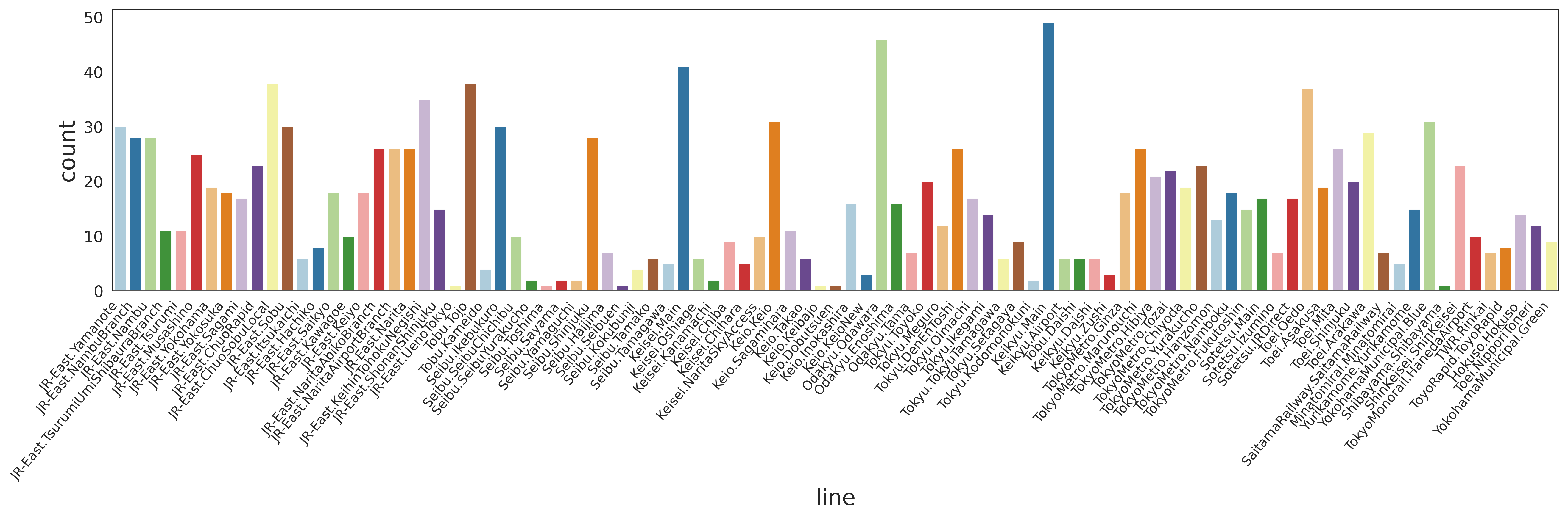}
\caption{The statistics of the station count of the 90 railway lines}
\label{fig:line90}
\end{figure*}
\subsection{Objective and Output}
Let the \textbf{matrix $\bm{K} \in \mathbb{R}^{n \times n}$ be the optimized traffic evacuation matrix} after the shutdown of train stations, where $k_{i,j}$ denotes the count of passengers that need to be evacuated from the $i$-th station to the $j$-th station. The objective is to find a traffic assignment output \( \textbf{K} = \{k_{ij}\} \) to evacuate the flow of passengers from affected stations to target stations. 
The solution must ensure efficient passenger evacuation while preserving network stability and minimizing evacuation cost. The evacuated strategy must satisfy the following constraints and optimization objectives.
\subsection{Constraints}
\noindent \textbf{(1) Equality Constraint:}The total number of evacuated passengers must match the affected passenger volume:  
\begin{equation}  
    \sum_{j} k_{ij} = y_{i}, j \in \hat{J},  
\end{equation}  
where \(y_{i} \) denotes the passenger volume at the affected station \( v_i \), and \( \hat{J} \) represents the candidate evacuation stations \( v_j \) that are geographically and topologically nearest to \( v_i \) and can accept evacuated passengers.  

\begin{equation}
    \bm{\hat{J}} \odot \bm{Y} 
\;=\;
\bm{K} \,\bm{J}.
\end{equation}
This matrix computation ensures an equivalent relation between $\bm{K}$ and the given data, which means that all people in blocked stations must leave for other stations.

\noindent \textbf{(2) Inequality Constraint:}
\textbf{Capacity Constraint:} All the stations have a common constant ratio for upper limit capacity:  
\begin{equation}  
     y_i' \leq x_i y_i, 
\end{equation}  
where \( y_i \) is the original capacity of $i$-th station, \( y_i' \) is the expanded capacity which defined by $x_{i}$, denoted as the upper limitation ratio.
\begin{equation}
  \bm{Y} \odot (\bm{J} - \bm{\hat{J}})
\;+\;
\bm{K}^{\top} \, (\bm{J} - \bm{\hat{J}})
\;\;\le\;\;
\bm{X} \,\odot\, \bm{Y} \,\odot\, (\bm{J} - \bm{\hat{J}}).
\label{eq:7}
\end{equation}
Equation~\ref{eq:7} can be interpreted as imposing a limit on the maximum volume at each station. It means that each blocked station has a fixed capacity $x_i$ in $\bm{X}$ that determines the maximum number of evacuated passengers it can accommodate.

\noindent \textbf{(3) Diagonal Constraint:}
\begin{equation}
    k_{i,i} = 0,
\quad
\forall \, i = 1,\dots,n.
\end{equation}
Each diagonal entry $k_{i,i}$ of $\bm{K}$ must be zero, which means that travelers cannot evacuate from their original stations to the same stations.

\noindent \textbf{(4) Additional Row/Column Constraints Depending on \(\bm{\hat{J}}\):}

This rule provides a more detailed constraint based on prior knowledge. For each \(i \neq j\), the values of \(k_{i,j}\) and \(k_{j,i}\) must satisfy certain bounds that depend on whether \(\hat{j}_{i,0}\) is 0 or 1, defined as follows in Equation~\ref{eq:9}:
\begin{equation}
\begin{cases}
k_{j,i} \;\ge\; 0, \quad k_{i,j} = 0 &\quad  \text{if $\hat{j}_{i,0} = 0$},\\
k_{i,j} \;\ge\; 0, \quad k_{j,i} = 0 &\quad  \text{if $\hat{j}_{i,0}= 1$}.\\
\end{cases}
\label{eq:9}
\end{equation}
It denotes that, for station $i$ in the normal status: inflow is possible; outflow is not allowed.
otherwise, station $i$ is blocked, outflow is mandatory and inflow is not allowed.

\subsection{Optimization}
Based on the assumptions and constraints, we formulate the following nonlinear problem for optimization.

\textbf{Travel Cost Minimization:} To minimize the travel time or walking distance caused by evacuation, we have $\bm{A_{cost}}$ representing the travel cost matrix within all stations (each element denoted as  \( a_{ij} \)). The optimization seeks to minimize the average travel cost for each passenger that need to be evacuated:  
\begin{equation}  
    \min \sum_{i \in \hat{J}} a_{ij} \cdot k_{ij}.  
\end{equation}  

The optimization process is based on matrix computations. The objective function is defined as the mean value of the element-wise product between $\bm{A_{cost}}$ and $\bm{K}$ (denoted as $(\bm{A_{cost}} \odot \bm{K})$), representing the average travel cost per person according to the decision matrix $\bm{K}$. This product is then multiplied  by $\bm{J}$ and its transpose $\bm{J}^{\top}$ to obtain a scalar value:
\begin{equation}
sum_{AK} = \bm{J}^{\top}(\bm{A_{cost}} \odot \bm{K}) \,\bm{J}.
\end{equation}
$(\cdot)^{\top}$ denotes matrix transpose. Similarly, we can easily represent the total number of evacuated passenger flow:
\begin{equation}
sum_{K} = \bm{J}^{\top} \bm{K} \,\bm{J}.
\end{equation}

The objective function is then given by
\begin{equation}
obj(\bm{K}) = \frac{sum_{AK}}{sum_{K} + \varepsilon}
\;=\;
\frac{\bm{J}^{\top}(\bm{A_{cost}} \odot \bm{K}) \,\bm{J}}{\bm{J}^{\top}\bm{K} \,\bm{J} + \varepsilon}.
\end{equation}
A small positive constant $\varepsilon$ was used to prevent division by zero (set to be $1 \times 10^{-6}$). 

The problem can be formulated as a multi-objective optimization problem that simultaneously minimizes both the \textbf{average evacuation overload} and the \textbf{total travel cost}, subject to flow conservation and capacity constraints, defined as the following equation~\ref{eq:14}.\begin{equation}
\begin{aligned}
&\min_{\bm{K}}\; obj(\bm{K}),\\
& \text{subject to Constraints (1), (2), (3), (4)}.
\end{aligned}
\label{eq:14}
\end{equation}

A suitable solver (e.g., an interior-point method) can be used to solve the above optimization problem and return the optimal matrix \(\bm{K}\).


\section{Experiment}
The data preprocessing was developed using Python 3.8.8, Numpy1.21.0, networkx 2.5, etc. The nonlinear optimization was developed using CasADi 3.6.7\footnote{https://web.casadi.org/}. The program was deployed on a Linux server (Architecture: x86\_64; CPU: 128 AMD Ryzen Threadripper
3990X 64-Core Processor).


\subsection{Study Area and Datasets}
\begin{table}[h]
\centering
\caption{Statistics of the dataset}
\begin{tabular}{cp{1.0cm}p{1.0cm}p{1.0cm}p{1.0cm}p{1.0cm}p{1.0cm}}
\toprule
Dataset   & \#Nodes & \#Lines& \#Operators\\\hline
Great Tokyo Railway Network  & 1113 & 90& 21\\
\bottomrule
\end{tabular}
\label{tab:data_sta}
\end{table}

We conducted the case study on the the railway network of Greater Tokyo, which is dominated by the world's most extensive urban rail network. The Greater Tokyo area is made up of Tokyo and the three neighboring prefectures of Saitama, Chiba, and Kanagawa. This area is home to around 30\% of Japan's total population. It is the most populated metropolitan area in the world, with a population of 41 million in 2024 \cite{WikipediaTokyo}. The network has a high rail transit usage ratio, with up to 48\% of the residents preferring to travel by train or metro services \cite{Jeph2022}. The Greater Tokyo's railway system, demonstrates high intensity and demands greater resilience. In this research, we selected central area of the railway network. The dataset statistics are summarized in Table~\ref{tab:data_sta}. Figure~\ref{fig:line90} gives an overview of station count on 90 lines. We collect static data (e.g. station name, line name, latitude and longitude) of the railway network from the website ekidata.jp\footnote{https://ekidata.jp/}. We also used passenger count data for the year 2019 \cite{hangli2022multi,ge2024k}. The traffic data are the average count of passengers per day in one year\footnote{https://www.toukei.metro.tokyo.lg.jp/homepage/ENGLISH.htm}.
\subsection{Experimental Setting}


In the simulation, $T_{lm}$ is set to 30 minutes, which means the disruption duration is 30 minutes. The affected passenger count is calculated by dividing the daily average by 20 hours*2, where 20 represents the daily railway operation time (in hours). The model represents the number of passengers affected by the emergencies on the rail system who need to be evacuated. In the experiment, $X$ represents the ratio of the maximum passenger volume capacity, and it was set to 1.5.

\section{Evaluation}


\subsection{Evaluation Metrics}

We evaluate the simulation performances using three following metrics:
\begin{itemize}
    \item EPF (Evacuated Passenger Flow) is count of evacuated passengers which was optimized from the origin station to target stations.
    \item PTT (Passenger Travel (Cost) Time) is of the passengers' travel time which was originally modeled in cost matrix from the origin station to all other target stations.
    \item ATT (Average Travel (Cost) Time) is calculated as the total sum of cost time of origin station, to be divided by the total sum of evacuated passenger count.
\end{itemize}

\subsection{Results}
\begin{figure}[h]
    \centering
    \includegraphics[height=1.6in, width=0.9\linewidth]{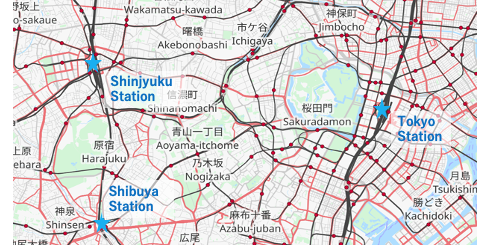} %
    \caption{Three stations were highlighted as {\LARGE\color{cyan} $\star$}; The map is based on the transportation view of OpenStreetMap, where red dots show  stations, and the black lines indicate railway lines.}
    \label{osm}
\end{figure}

We conducted simulations on three major stations: Tokyo, Shinjuku, and Shibuya (Fig.~\ref{osm}). In each simulation, we assume that one of the corresponding station is disrupted for $T_{lm}$ (30 minutes), and all passengers from that station must be evacuated to other stations.  

\begin{table}[h]
\centering
\caption{ATTs of Tokyo, Shinjuku, and Shibuya}
\begin{tabular}{p{3.0cm}p{1.0cm}p{1.0cm}p{1.0cm}}
\toprule
Station   & Tokyo &Shinjuku  &Shibuya \\\hline
ATT (mins)   & 4.2& 8.3& 7.6 \\
\bottomrule
\end{tabular}
\label{tab:results}
\end{table}
Extracted from the statistic data, the effected passenger count within $T_{lm}$ of Tokyo, Shinjuku, and Shibuya are 33488, 81639 and 67609, respectively. The ATT (listed in Table~\ref{tab:results}) demonstrates the final average travel time can be optimized to under 10 minutes, meaning it takes on average less than 10 minutes to reassign passengers for evacuation. 
\begin{figure}[h]
    \centering
    \includegraphics[height=1.2in, width=0.95\linewidth]{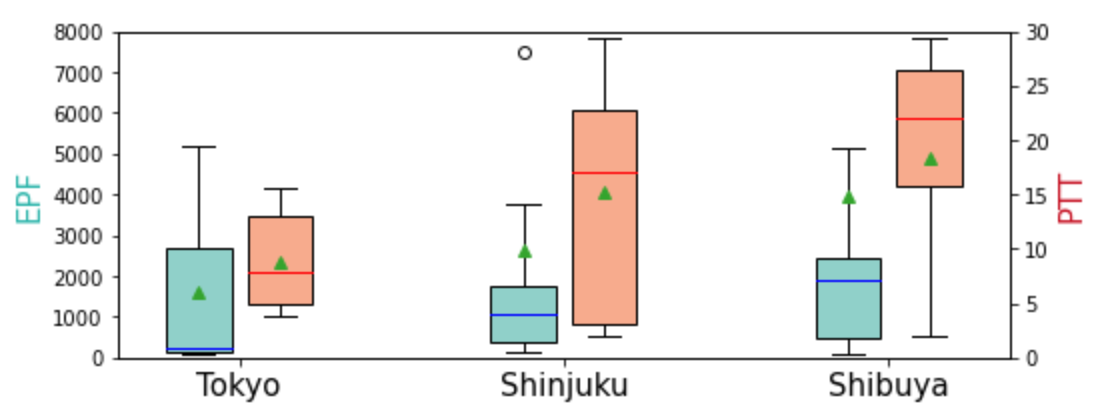} 
    \caption{Box plots of station-level re-assigned passenger count (EPF) and travel time (PTT).}
    \label{box}
\end{figure}

In Fig.~\ref{box}, we can see that the median number of evacuated passengers (EPF) at the three stations is 243, 1070, and 1883, respectively (marked as blue line). The average values are 1629, 2631, 3977, respectively (marked by green triangle). On the other hand, the average passenger travel times (PTT) are 7.8, 17.1, and 22.1 minutes. In summary, Tokyo station is the easiest station to evacuate passengers from due to the abundance of nearby target stations, as observed in  Fig.~\ref{osm}. In contrast, passengers at Shinjuku station must travel longer distances to reach alternative stations. Meanwhile, passengers evacuated from Shibuya station travel the longest distances, as there are fewer surrounding stations.

\begin{figure}[htbp]
    \centering
    \subfloat[Evacuated from Tokyo]{\includegraphics[height=1.3in, width=1\linewidth]{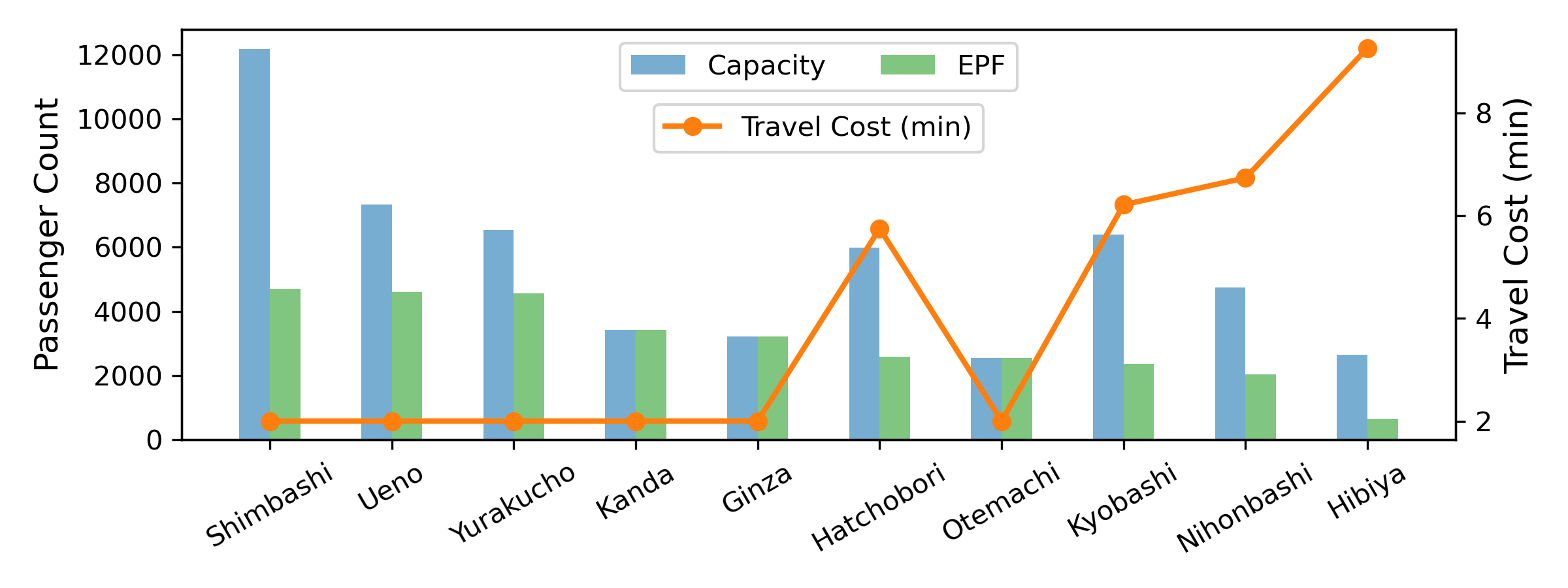}\label{tokyo}} \\
    \subfloat[Evacuated from Shinjuku]{\includegraphics[height=1.3in, width=1\linewidth]{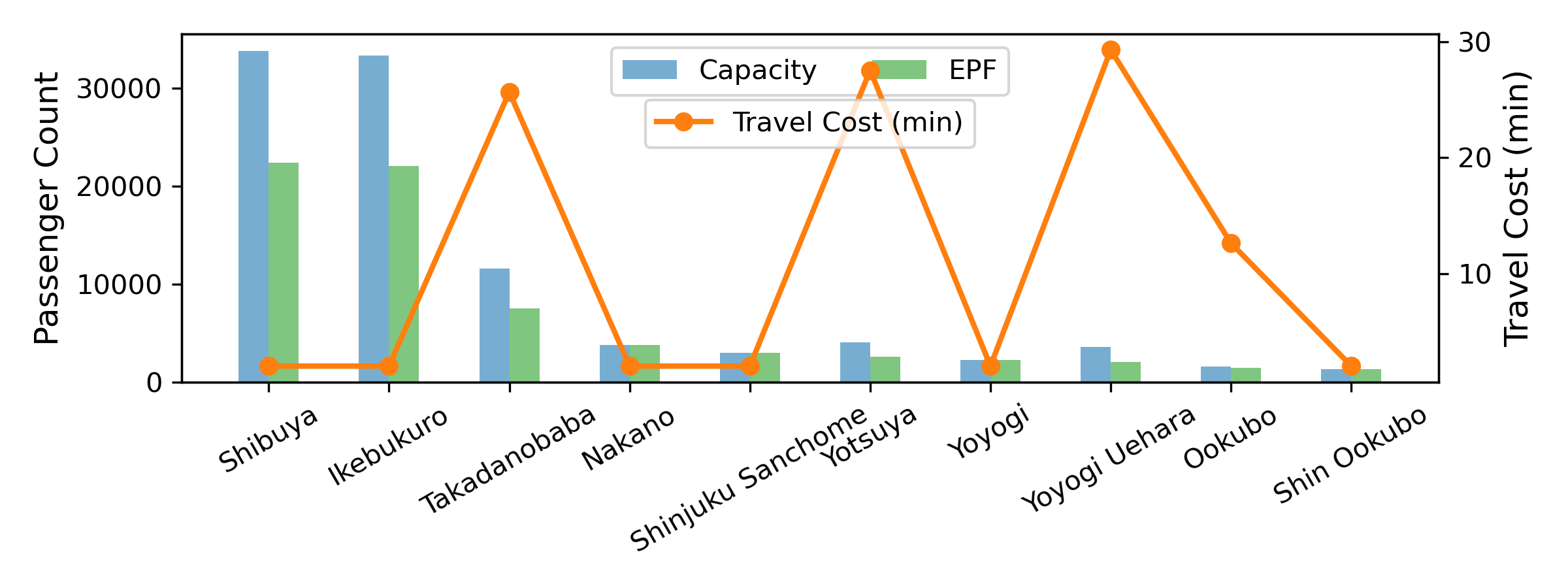}\label{shinjyuku}} \\
    \subfloat[Evacuated from Shibuya]{\includegraphics[height=1.3in, width=1\linewidth]{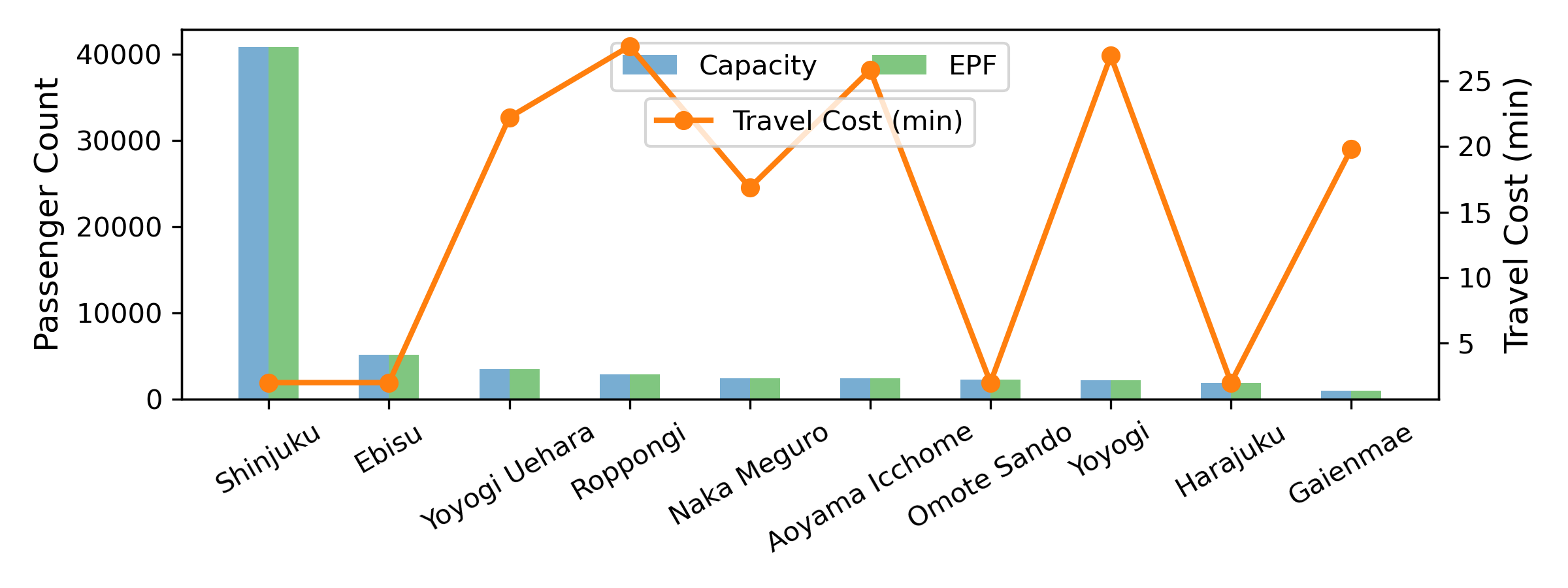}\label{shibuya}} 
    \caption{Capacity and EPF of top 10 stations from the origins of Tokyo station, Shinjuku station and Shibuya station), along with the corresponding travel cost.}
    \label{stations}
\end{figure}

We visualized the maximum capacity of station and evacuated passenger counts and travel times to each station. In Fig.~\ref{stations}, we can observe that directly linked stations are always the preferred choices, as  passengers can travel to these stations in a very short time (marked the travel cost as 2 minutes as we set up in our cost matrix before). Compared to the other two stations, Shinjuku and Shibuya, Tokyo (shown in Fig.~\ref{tokyo}) station has more connectivities makes the majority of the evacuated target stations are the connected stations, followed by the geographically closest stations. In the case of Shinjuku (Fig.~\ref{shinjyuku}), some relatively large capacity but unconnected stations are chosen as evacuation target. Furthermore, in Fig.~\ref{shibuya}, we can oberve that the first selected station (i.e. Shinjuku) receives a large number of passengers up to station capacity, while the others only receive passenger counts in the thousands. This is because Shinjuku, being a very large neighboring station, absorbs a significant share of passengers.

\section{Conclusion}
In this study, we present a traffic evacuation model for railway disruptions aimed at improving resilience. Our model focuses on optimizing evacuation strategies in response to disruptions affecting multiple nodes within the rail network. We define a multi-weighted, geographically feature-fused travel cost matrix that dynamically adjusts based on predefined disruption durations. By setting appropriate thresholds, we can quickly identify candidate stations for evacuation. We developed an optimization model using nonlinear programming that considers both the capacity of surrounding evacuation stations and the associated travel costs. In addition, our model can flexibly simulate disruptions across any number of stations and lines, enhancing its applicability. We simulated our model using Greater Tokyo's large-scale railway network and tested three major stations. The results show that evacuation times during emergencies are optimized and remain within acceptable recovery periods, demonstrating the effectiveness of our approach. 



\section*{Acknowledgment}
The work of Francesco Flammini was partly supported by the Swiss State Secretariat for Education, Research and Innovation (SERI) under contract no. 24.00528 (PhDs EU-Rail project). The project has been selected within the European Union’s Horizon Europe research and innovation programme under grant agreement no. 101175856. Views and opinions expressed are however those of the authors only and do not necessarily reflect those of the funding agencies, which cannot be held responsible for them.

\bibliography{IEEEabrv,ref}

\bibliographystyle{IEEEtran}

\end{document}